\documentstyle[12pt,amssymb,amstex]{amsart}
\numberwithin{equation}{section}

\def\ca{{\cal A}}
\def\cb{{\cal B}}

\def\ch{{\cal H}}

\def\cl{{\cal L}}
\def\cam{{\cal M}}
\def\cn{{\cal N}}

\def\car{{\cal R}}

\def\ga{{\frak A}}

\def\gam{{\frak M}}
\def\gn{{\frak N}}

\def\gar{{\frak R}}

\def\gz{{\frak Z}}

\def\bc{{\Bbb C}}

\def\a{\alpha}
\def\b{\beta}
  \def\G{\Gamma}
  \def\D{\Delta}
\def\eps{\epsilon}
\def\eeps{\varepsilon}
 
\def\m{\mu}

\def\n{\nu}

\def\s{\sigma}
\def\t{\tau}
\def\ff{\varphi}\def\f{\phi} 
\def\Th{\Theta}
 \def\Om{\Omega}

\def\dim{\mathop{\rm dim}}
\def\ots{\overline{\otimes}}

\def\nn{\nonumber}
\def\di{\mathop{\rm d}\!}

\newcommand{\ty}[1]{\mathop{\rm {#1}}}

\newtheorem{thm}{Theorem}
\newtheorem{lem}{Lemma}
\newtheorem{cor}{Corollary}
\newtheorem{prop}{Proposition}

\newtheorem{defin}{Definition}

\begin{document}

\title[tensor products over subalgebras]
{the structure of the $W^{*}$--tensor product over a 
$W^{*}$--subalgebra and its predual 
($\s$--finite case)} 
\author{Francesco Fidaleo}
\address{Francesco Fidaleo\\
Dipartimento di Matematica\\
Universit\`{a} di Roma ``Tor Vergata''\\
Via della Ricerca Scientifica 1, 00133 Roma, Italy.}
\email{{\tt fidaleo@@axp.mat.uniroma2.it}}
\begin{abstract}
Let $M$, $N$, $R$ be $W^{*}$--algebras, with $R$ unitally embedded 
in both $M$ and $N$. by using Reduction Theory, we extend the previous description of the
$W^{*}$--tensor product 
$M\ots_{R}N$ over the common $W^{*}$--subalgebra $R$ and its 
predual $\big(M\ots_{R}N\big)_{*}$ to the $\s$--finite case.
\end{abstract}

\maketitle

\section{introduction}
 
The structure of 
$W^{*}$--tensor products over common subalgebras 
was investigated in \cite{SZ}. The tensor product over 
subalgebras generalizes the usual  
$W^{*}$--tensor product $M\ots N$, corresponding to the case when $R$ 
is the field of complex numbers. 
The main goal of the paper \cite{SZ} is to prove the 
generalization of the celebrated Tomita Commutation Theorem to 
the case of tensor products over subalgebras in the full generality. 
If properly interpreted, the commutant assumes the form
$$
\big(M\ots_{R}N\big)'=M'\ots_{R'}N'\, ,
$$    
and reduces itself to Tomita Commutation Theorem when $R=\bc$ acting 
on $\bc$, see \cite{SZ}, Section 5.9.  

In the case of $W^{*}$--tensor product $M\ots N$, there is an explicit 
description of the predual as
\begin{equation}
\label{41}
\big(M\ots N\big)_{*}\cong M_{*}\widehat\otimes N_{*}
\end{equation}
where $\widehat\otimes$ denotes the operator--projective tensor 
product, see \cite{ER1}.\footnote{For standard definitions and 
results relative to operator space theory, the reader is referred to 
\cite{E} and the references cited therein.}
Notice that   
\eqref{41} has several applications to various fields, 
see e.g. \cite{BDL, F1, F, Su} and the reference cited therein.

In the present paper we extend to the $\s$--finite case, 
the results of the previous paper 
\cite{F2} relative to the structure of the $W^{*}$--tensor product 
of von Neumann algebras over a common $W^{*}$--subalgebra. This is 
done by using reduction theory. Namely, 
using the extremal decomposition of KMS states (see e.g. \cite{A}) 
for $C^{*}$--dynamical systems, and the general properties of 
$W^{*}$--tensor products over common subalgebras, we describe the 
structure of such a tensor product, following the lines adopted in 
\cite{R}. Also for this situation, $R$, $M$, $N$, and finally 
$M\ots_{R}N$ admit a common decomposition over the centre of $R$. 
Theorem 1 and Corollary 1 of \cite{F2} can be generalized to $\s$--finite 
case, and describe the structure of the $W^{*}$--tensor product
a common $W^{*}$--subalgebra, and its predual respectively. These 
results can 
be considered as a further step towards the fully understand of the 
general case.

\section{on the tensor product over subalgebras}

We start with 
the definition of the tensor product over a subalgebra, as well as 
some preliminary properties already considered in \cite{SZ}.
\begin{defin}
\label{51}
Let $I_{A}\in R\subset M\subset A$, $I_{A}\in R\subset N\subset A$ 
be inclusions of  von Neumann algebras. The $W^{*}$--algebra $A$ 
is said to be the $W^{*}$--tensor product of $M$ and $N$ over their common 
$W^{*}$--subalgebra $R$ if
\begin{itemize}
\item[(i)] $A=M\bigvee N$,
\item[(ii)] for some normal 
representation $\pi$ of $A$, 
$$
\pi(M)\subset\pi(R)\bigvee M_{1}\, ,\quad 
\pi(N)\subset\pi(R)\bigvee N_{1}
$$ 
for commuting type $\ty{I}$ von Neumann subalgebras 
$M_{1},N_{1}\subset\pi(R)'$ whose common centre $Z$ coincides with that 
of $R$.
\end{itemize}
\end{defin}

In the situation described by Definition \ref{51}, we write 
$A=M\ots_{R}N$ and identify $\pi(A)$ with $A$ itself. Any such a 
representation as that in Definition \ref{51} is said to be a 
{\it splitting} one.
Further, the type $\ty{I}$ von Neumann algebras
$M_{1}$, $N_{1}$ given in (ii) can be chosen to be homogeneous, see 
\cite{F2}, Proposition 1. Then, we can start without loss of 
generality, by the following inclusions of von Neumann algebras 
\begin{align}
\label{1}
&I_{1}\otimes R\subset M\subset B(H_{1})\ots R\, ,\nn\\
&R\otimes I_{2}\subset N\subset R\ots B(H_{2})\, ,
\end{align}
where $R$ is acting on $B(H)$.
The $W^{*}$--tensor product of $M$ and $N$ 
over $R$ is described with a slightly abuse of 
notation, by 
\begin{equation}
\label{1a1}
M\ots_{R}N=\big(M\otimes I_{2}\big)\bigvee \big(I_{1}\otimes N\big)\, .
\end{equation}

Further,
$$
I_{1}\otimes Z(R)\subset Z(M)\, ,\quad Z(R)\otimes I_{2}\subset Z(N)\, .
$$

From now on, we suppose that $A$ is $\s$--finite if it is not otherwise 
specified.

We construct normal faithful conditional expectations which are useful 
in the sequel.\footnote{We sometimes identify $M$, $N$, $R$ with their 
isomorphic copies $M\otimes I_{2}$, $I_{1}\otimes N$ and 
$I_{1}\otimes R \otimes I_{2}$ in $A\equiv M\ots_{R}N$.}
\begin{prop}
\label{coex}
There are normal faithful conditional expectations
$$
\eeps_{1}:A\mapsto M\,,\quad \eeps_{2}:A\mapsto N\,.
$$

Furthermore, 
\begin{equation}
\label{1cex}
\eeps_{1}\circ\eeps_{2}=\eeps_{2}\circ\eeps_{1}.
\end{equation}
\end{prop}
\begin{pf}
Consider for $\ff\in L^{1}(H_{2})_{+,1}$, the slice map of 
$B(H_{1})\ots R\ots B(H_{2})$ onto $B(H_{1})\ots R$,
$$
F^{1}_{\ff}(x\otimes y)=\ff(y)x\otimes I_{2}\,,
$$
together with its restriction $E^{1}_{\ff}:=F^{1}_{\ff}\lceil_{A}$. The 
set $\{E^{1}_{\ff}\,|\,\ff\in L^{1}(H_{2})_{+,1}\}$ is a separating 
family of conditional expectations of $A$ onto $M$. As $A$ is 
$\s$--finite, there exists a denumerable maximal family 
$\{E^{1}_{\ff_{k}}\}$ with mutually orthogonal support--projections (see 
\cite{St}, 11.5). Define 
\begin{equation*}
\eeps_{1}:=\sum_{k=1}^{\infty}\frac{1}{2^{k}}E^{1}_{\ff_{k}}\,,
\end{equation*}
which is the searched conditional expectation. Starting from
$$
F^{2}_{\psi}(x\otimes y)=\psi(x)I_{1}\otimes y\,,
$$
construct 
$\eeps_{2}:A\mapsto N$ 
by a denumerable maximal family 
$\{E^{2}_{\psi_{k}}\}$ with mutually orthogonal support--projections, 
where $\psi_{k}\in L^{1}(H_{1})_{+,1}$.
We have, for $a\in A$,
\begin{align*}
&E^{1}_{\ff}(E^{2}_{\psi}(a))\equiv F^{1}_{\ff}(F^{2}_{\psi}(a))\\
=&F^{2}_{\psi}(F^{1}_{\ff}(a))\equiv E^{2}_{\psi}(E^{1}_{\ff}(a))\,,
\end{align*}
which leads to the assertion taking into account \eqref{1cex} for 
$\eeps_{1}$, and the analogous one for $\eps_{2}$.
\end{pf}

Define 
\begin{equation}
\label{2cex}
\eps:=\eeps_{1}\circ\eeps_{2}\,,\quad\eps_{1}:=\eeps_{1}\lceil_{M}\,,\quad
\eps_{2}:=\eeps_{2}\lceil_{N}\,.
\end{equation}

Consider the diagram
\newpage
\begin{figure}[hbt]
\begin{center}
\begin{equation}
\label{edia}
\begin{picture}(350,120)
\put(100,20){\begin{picture}(130,100)
\thicklines
\put(42,83){$\vector(0,-1){48}$}
\put(122,83){$\vector(0,-1){48}$}
\put(110,23){$\vector(-1,0){55}$}
\put(105,80){$\vector(-1,-1){43}$}
\put(110,94){$\vector(-1,0){55}$}
\thinlines
\put(35,90){$M$}
\put(35,20){$R$}
\put(118,20){$N$}
\put(118,90){$M\ots_{R}N$}
\put(27,60){$\eps_{1}$}
\put(128,60){$\eeps_{2}$}
\put(80,12){$\eps_{2}$}
\put(80,100){$\eeps_{1}$}
\put(80,65){$\eps$}
\end{picture}}
\end{picture}
\end{equation}
\end{center}
\end{figure}
\noindent
where $\eeps_{j}$, $j=1,2$ are given in Proposition \ref{coex}, 
and $\eps$, $\eeps_{j}$, $j=1,2$ are given in \eqref{2cex}. 
\begin{cor}
The above diagram  
gives rise to a commuting square of normal faithful conditional 
expectations.\footnote{For the 
definition of a commuting square of conditional expectations see 
e.g. \cite{J}.}
\end{cor}
\begin{pf}
The proof immediately follows by \eqref{1cex} as 
$M\otimes I_{2}\bigwedge I_{1}\otimes N=I_{1}\otimes R \otimes I_{2}$.
\end{pf}

Now we pass to investigate the standard representation of 
$M\ots_{R}N$ for general (non necessarily $\s$--finite) 
$W^{*}$--algebras.

By applying the considerations in the beginning of Section 6 
of \cite{SZ}, we can describe $M\ots_{R}N$ in the following way.
Put 
\begin{equation*}
\tilde M:=M\bigwedge B(H_{1})\ots Z\,,\quad
\tilde N:=N\bigwedge Z\ots B(H_{2})\,. 
\end{equation*}

Then,
\begin{equation*}
M\ots_{R}N=\tilde M\otimes I_{2}\bigvee R\bigvee I_{1}\otimes\tilde N
\equiv\tilde M\ots_{Z}R\ots_{Z}\tilde N\,.
\end{equation*}
\begin{thm}
\label{sta}
The standard representation of the $W^{*}$--tensor product over a 
$W^{*}$--subalgebra is a splitting representation.
\end{thm}
\begin{pf} 
After taking a possible ampliation, we can suppose that $R$ is acting 
on $L^{2}(Z)\otimes H$,\footnote{In the general 
commutative case, representing $Z$ by the GNS representation 
relative to a normal faithful weight $\ff$ (\cite{SZ}, Section 10.14), we have
$Z\sim L^{\infty}(\G_{\ff},\m_{\ff})$ where $\m_{\ff}$ is a positive Radon measure 
on the locally compact dense 
subspace $\G_{\ff}$ of the spectrum $\Om$ of $M$, see \cite{T}, Theorem 
III.1.18. Hence, $L^{2}(Z)\cong L^{2}(\G_{\ff},\m_{\ff})$.}
and the standard representations of
$\tilde M$, $\tilde N$ and 
$R$ can be obtained by induction on the Hilbert spaces  
$H_{1}\otimes L^{2}(Z)$, $L^{2}(Z)\otimes H_{2}$, $L^{2}(Z)\otimes H$
where they are naturally acting. Let $e_{1}\in\tilde M'$, $e_{2}\in\tilde N'$, 
$e\in R'$ be the corresponding selfadjoint projections. Put 
$$
E:=(e_{1}\otimes I\otimes I_{2})(I_{1}\otimes e\otimes I_{2})
(I_{1}\otimes I\otimes e_{2})\,.
$$

Due to Theorem 4.7, the operator defined above is a selfadjoint 
projection. Furthermore, $e_{1}\in B(H_{1})\ots Z$, 
$e_{2}\in Z\ots B(H_{2})$. 

The proof follows as $(M\ots_{R}N)E$ is the 
standard representation of $M\ots_{R}N$ acting on 
$E(H_{1}\otimes L^{2}(Z)\otimes H\otimes H_{2})$, and the 
ampliations--inductions  
$(e_{1}(B(H_{1})\ots Z)e_{1}\otimes I\otimes I_{2})E$,
$(I_{1}\otimes I\otimes e_{2}(Z\ots B(H_{2}))e_{2})E$ of the
reduced algebras $e_{1}(B(H_{1})\ots Z)e_{1}$,
$e_{2}(Z\ots B(H_{2}))e_{2}$ are the splitting type $\ty{I}$ 
algebras appearing in Definition \ref{51}
\end{pf}

We end by noticing that the standard representation is not 
homogeneous in general.

\section{on the decomposition of von Neumann algebras arising from left 
Hilbert algebras and their predual}

Let $\{\ga_{\xi}\}_{\xi\in\Om}$ be a field of left Hilbert 
algebras defined on the finite measure space $(\Om,\m)$. We suppose 
that such a field of left Hilbert algebras satisfies Conditions
(1.1)--(1.5) listed in \cite{R}. In this situation, the field 
$\{\ga_{\xi}\}_{\xi\in\Om}$ is said to be {\it integrable}. 
This is the case arising from left Hilbert 
algebras associated to GNS representations of KMS states for 
$C^{*}$--dynamical systems.
To the field mentioned above, there is associated a Hilbert space
$\ch$ which is the direct integral of an integrable field 
$\{\ch_{\xi}\}_{\xi\in\Om}$ of Hilbert spaces, see e.g. \cite{W}.

Let $\cl(\ga)$ be the associated left von Neumann algebra
together with the corresponding field $\{\cl(\ga_{\xi})\}_{\xi\in\Om}$ of 
left von Neumann algebras. Consider a measurable vector field
$\xi\mapsto T(\xi)$ such that $T(\xi)\in\cl(\ga_{\xi})'$ almost 
surely. It is readily seen that 
\begin{equation}
\label{diop}
T:=\int^{\oplus}_{X}T(\xi)\di\m(\xi)
\end{equation}
defines an element of $\cl(\ga)'$. Conversely, any $T\in\cl(\ga)'$ 
admits a (essentially unique) decomposition as above. 
Indeed, denote by $\gz$ the algebra 
consisting of all diagonal operators on $\ch$ (see e.g. \cite{T, VW}). 
Any such $T$ is decomposable as $T\in\gz'$, see Theorem 1.7 of 
\cite{W}, and Proposition 2.1 of \cite{VW}. Define $\Th:=JTJ\in\cl(\ga)$, where 
$J$ is the canonical conjugation associated to $\ga$ (see e.g. 
\cite{SZ1}, Section 10.1). By Proposition 1.3 and Theorem 1.4 of \cite{R}, $\Th$ 
admits a unique natural decomposition 
$\xi\mapsto\Th(\xi)$.\footnote{A natural decomposition for $T\in\cl(\ga)$
(resp $T\in\cl(\ga)'$) is a decomposition $\xi\mapsto T(\xi)$ of $T$
such that $T(\xi)\in\cl(\ga_{\xi})$ (resp $T(\xi)\in\cl(\ga_{\xi})'$)
almost surely, see \cite{R}, Definition 1.2.} Taking into account 
(1.4) of \cite{R}, $\xi\mapsto J(\xi)\Th(\xi)J(\xi)$ provides a 
natural decomposition of $T$ as in \eqref{diop}, which is essentially unique.  

Now we pass to the description of the predual $\cl(\ga)_{*}$ of 
$\cl(\ga)$ in the situation described above. This description 
parallels the analogous one concerning the separable situation (see 
e.g. \cite{T}, Section IV. 8), taking into account of appropriate changes.

Consider the subfield of ${\displaystyle\prod_{\xi\in X}\cl(\ga_{\xi})_{*}}$
consisting of elements $\xi\mapsto\ff(\xi)$ such that
\begin{itemize}
\item[(i)] the map $\xi\mapsto\ff(\xi)(T(\xi))$ is measurable for each 
$T\in\cl(\ga)$, ${\displaystyle T=\int^{\oplus}_{X}T(\xi)\di\m(\xi)}$ 
being the natural decomposition of $T$; 
\item[(ii)] there exists an element $c_{\ff}\in L^{1}(X,\m)_{+}$ such 
that $\|\ff(\xi)\|\leq c_{\ff}(\xi)$ almost 
everywhere.\footnote{Notice that we cannot conclude, 
in non separable cases, that $\xi\mapsto\|\ff(\xi)\|$ is measurable.}
\end{itemize}
\begin{prop}
There is a one--to--one correspondence between elements $\ff\in\cl(\ga)_{*}$ 
and elements $\xi\in X\mapsto\ff(\xi)\in\cl(\ga_{\xi})_{*}$ satisfying 
$\mathop{(\rm{i}), (\rm{ii})}$ above.
\end{prop} 
\begin{pf}
Let $\xi\mapsto\ff(\xi)$ be a measurable field of functional as 
above. Define
$$
\ff(T):=\int_{X}\ff(\xi)(T(\xi))\di\m(\xi)\,,
$$
which is well defined by Proposition 1.3 of \cite{R}. 
We get
$$
|\ff(\xi)(T(\xi))|\leq\|\ff(\xi)\|\|T\|\leq c_{\ff}(\xi)\|T\|\,,
$$
which means
$$
|\ff(T)|\leq\int_{X}|\ff(\xi)(T(\xi))|\di\m(\xi)
\leq\bigg(\int_{X}c_{\ff}(\xi)\di\m(\xi)\bigg)\|T\|\,,
$$
that is $\|\ff\|\leq\|c_{\ff}\|$. It is readily seen by 
Dominated Convergence Theorem, that $\ff$ is normal. Moreover, by 
considering the polar decomposition of normal functionals (see e.g. 
Theorem 5.16 of \cite{SZ1}), and Theorem 1.4 of \cite{R}, if $\ff$ is 
the null functional, then $\ff(\xi)=0$ almost surely. The construction 
of a field of functional $\xi\mapsto\ff(\xi)$ as above, starting from  
$\ff\in\cl(\ga)_{*}$ follows the same line of the analogous 
construction in the separable situation, see \cite{T}, Proposition IV. 
8. 34.
\end{pf}

Summarizing, we have the following terminology. Define
\begin{align*}
&M:=\cl(\ga)\,,\qquad M(\xi):=\cl(\ga_{\xi})\,,\,\,\,\quad\xi\in X\,\\
&M':=\cl(\ga)'\,,\,\,\,\quad M(\xi)':=\cl(\ga_{\xi})'\,,\,\quad\xi\in X\,\\
&M_{*}:=\cl(\ga)_{*}\,,\,\,\,\,\,\,\, 
M(\xi)_{*}:=\cl(\ga_{\xi})_{*}\,,\quad\xi\in X\,.
\end{align*}

We write
\begin{align*}
&M=\int^{\oplus}_{X}M(\xi)\di\m(\xi)\,,\\
&M'=\int^{\oplus}_{X}M(\xi)'\di\m(\xi)\,,\\
&M_{*}=\int^{\oplus}_{X}M(\xi)_{*}\di\m(\xi)\,.
\end{align*}

\section{the structure of the tensor product over a subalgebra and 
its predual}

We proved in Theorem \ref{sta}, that the standard representation is a 
splitting one. In general, it is non homogeneous. We start by 
recalling the structure of such a standard representation of
$A\equiv M\ots_{R}N$.
 
Let $I$, $J$ be the sets of 
cardinalities appearing in the decomposition in homogeneous parts 
(\cite{T}, Theorem V.1.27) of $M_{1}$, $N_{1}$ respectively appearing 
in the standard representation of $A$. We 
have
$$
H=\bigoplus_{\a,\b}H_{\a}\otimes H_{\a,\b}\otimes K_{\b}
$$
where $\dim(H_{\a})=\a$, $\dim(K_{\b})=\b$. Accordingly,
\begin{align}
\label{homo}
R=\bigoplus_{\a,\b}I_{H_{\a}}\otimes R_{\a,\b}\otimes 
I_{K_{\b}}\, , \quad
&Z=\bigoplus_{\a,\b}I_{H_{\a}}\otimes Z_{\a,\b}\otimes I_{K_{\b}}\,,\nn\\
M=\bigoplus_{\a,\b}M_{\a,\b}\otimes I_{K_{\b}}\, ,\quad
&N=\bigoplus_{\a,\b}I_{H_{\a}}\otimes N_{\a,\b}\, ,\\
M_{1}=\bigoplus_{\a,\b}\cb(H_{\a})\ots Z_{\a,\b}\otimes I_{K_{\b}}\, ,
\quad
&N_{1}=\bigoplus_{\a,\b}I_{H_{\a}}\otimes 
Z_{\a,\b}\ots\cb(K_{\b})\,.\nn
\end{align}

Looking at any single summand, we have for every $\a\in I$, $\b\in J$,
\begin{align*}
&I_{H_{\a}}\otimes R_{\a,\b}\subset M_{\a,\b}\subset
\cb(H_{\a})\ots R_{\a,\b}\, ,\\
&R_{\a,\b}\otimes I_{K_{\b}}\subset N_{\a,\b}\subset 
R_{\a,\b}\ots\cb(K_{\b})\,.
\end{align*}

As $A$ is $\s$--finite, we have that $I$ and $J$ are finite or 
countable sets. In this situation
\begin{equation}
\label{sttp}
M\ots_{R}N=\bigoplus_{\a,\b}\big(M_{\a,\b}\otimes I_{K_{\b}}\big)
\bigvee\big(I_{H_{\a}}\otimes N_{\a,\b}\big)\,.
\end{equation}

Let now $N\subset M$ be an inclusion of $\s$--finite von Neumann 
algebras such that there exists a normal faithful conditional 
expectation $E:M\mapsto N$ of $M$ onto $N$. Pick a normal faithful 
state $\ff$ on $N$, and consider $\psi:=\ff\circ E$. Consider the 
weakly dense $C^{*}$--subalgebra $\cam$ (resp. $\cn$)
made of elements $T\in M$
(resp. $T\in N$)  such that 
the function $\t\mapsto\s^{\psi}_{t}(T)$ (resp. 
$\t\mapsto\s^{\ff}_{t}(T)$) is continuous w.r.t. the norm 
topology.
\begin{lem}
\label{resex}
We have $\cn=\cam\bigcap N$ and $E\lceil_{\cam}$ is a conditional 
expectation of $\cam$ onto $\cn$.
\end{lem}
\begin{pf}
By Takesaki Theorem (see e.g. \cite{St}, Section 10), it is enough to 
show that $E(\cam)\subset\cn$. Let $M$ be directly represented in 
standard form on $\ch$, such that the standard vector $\Psi\in\ch$ 
gives rise the state $\psi$ on $M$. Let $\D$ be the modular operator 
relative to $\Psi$, and $P\in N'$ the projection inducing the standard 
representation of $N$ on $P\ch$. Then $P$ reduces $\D$, and we have
\begin{align*}
\s^{\ff}_{t}(E(T))\Psi=&\s^{\psi}_{t}(E(T))\Om=\D^{it}PT\Psi\\
=&P\D^{it}T\Psi=E(\s^{\psi}_{t}(T))\Psi\,,
\end{align*}
that is $\s^{\ff}_{t}\circ E=E\circ\s^{\psi}_{t}$. Namely, if $T$ is 
a regular element of $M$, $E(T)$ is a regular element of $N$.
\end{pf}

Fix a normal faithful state $\f$ on $R$, and extend it to all of $A$ 
by using the conditional expectation $\eps$. Taking into account the 
commuting square given in \eqref{edia}, we easily have
$$
\f\circ\eps\lceil_{M}=\f\circ\eps_{1}\,,\quad
\f\circ\eps\lceil_{N}=\f\circ\eps_{2}\,.
$$

Moreover, 
\begin{align}
\label{11a1}
[\ga_{\f\circ\eps}]=&[\gam_{\f\circ\eps_{1}}\bigcup\gn_{\f\circ\eps_{2}}]\nn\\
[\gar_{\f}]=&[\gam_{\f\circ\eps_{1}}\bigcap\gn_{\f\circ\eps_{2}}]\,,
\end{align}
where $\ga_{\f\circ\eps}$, $\gam_{\f\circ\eps_{1}}$, 
$\gn_{\f\circ\eps_{2}}$, $\gar_{\f}$ are the left Hilbert 
algebras with (the same) unity relative to the states $\f\circ\eps$,
$\f\circ\eps_{1}$, $\f\circ\eps_{2}$, $\f$, respectively, 
$[\,\cdot\,]$ denotes the closed generated subspace, and 
$\ga_{\f\circ\eps}$ is a dense subspace of the Hilbert space of the standard 
representation of $A$.

Consider the $C^{*}$--subalgebras $\ca$, $\cam$, $\cn$, $\car$ of 
regular elements of $A$, $M$, $N$, $R$ w.r.t. the modular group.
\begin{prop}
\label{regal}
We have
$$
M\ots_{R}N\equiv\pi_{\f\circ\eps}(\ca)''=\big(\pi_{\f\circ\eps_{1}}(\cam)''\big)
\ots_{\big(\pi_{\f}(\car)''\big)}\big(\pi_{\f\circ\eps_{2}}(\cn)''\big)\,.
$$
\end{prop}
\begin{pf}
The assertion immediately follows by Lemma \ref{resex}, and the 
previous considerations.
\end{pf}
  
Now we are ready to describe the structure of the $W^{*}$--tensor 
product over a $W^{*}$--subalgebra and its predual. In order to give a 
more readable description, we treat each homogeneous component 
separately.\footnote{We pursue such a strategy in order to give a more 
readable description of the structure of $M\ots_{R}N$ and its predual, 
see below.} 

Fix a sequence $\{\f_{\a\b}\}$ of normal states on $R_{\a\b}$, 
one for each homogeneous component in \eqref{homo}. By restricting 
ourselves to each homogeneous component, we consider separately 
$A_{\a\b}:=M_{\a\b}\ots_{R_{\a\b}}N_{\a\b}$. In this situation,
$A_{\a\b}$ has the form \eqref{1a1}, for inclusion of 
algebras as in \eqref{1}. Let $\f_{\a\b}\circ\eps$, 
$\f_{\a\b}\circ\eps_{1}$ and $\f_{\a\b}\circ\eps_{2}$ be the corresponding 
extensions to $A_{\a\b}$, $M_{\a\b}$ and $N_{\a\b}$ respectively. 
Consider the $C^{*}$--dynamical systems
$(\ca_{\a\b},\a^{\f_{\a\b}\circ\eps},\f_{\a\b}\circ\eps)$, 
$(\cam_{\a\b},\a^{\f_{\a\b}\circ\eps_{1}},\f_{\a\b}\circ\eps_{1})$
$(\cn_{\a\b},\a^{\f_{\a\b}\circ\eps_{2}},\f_{\a\b}\circ\eps_{2})$, 
$(\car_{\a\b},\a^{\f_{\a\b}},\f_{\a\b})$, 
where the $\a^{\#}$ are the restriction of the modular groups $\s^{\#}$ 
to the corresponding regular elements. Consider the extremal decomposition of 
the $\a^{\f_{\a\b}}$--KMS state $\f_{\a\b}$. 
Let $\m_{\a\b}$ be the maximal measure on the 
compact convex set $K(\car_{\a\b},\a^{\f_{\a\b}})$ corresponding to 
$\f_{\a\b}$, $K(\car_{\a\b},\a^{\f_{\a\b}})$ denoting all 
$\a^{\f_{\a\b}}$--KMS states on $\car_{\a\b}$. It is well known that 
$\m_{\a\b}$ coincides with the central measure of $\f_{\a\b}$ on 
$K(\car_{\a\b},\a^{\f_{\a\b}})$, see e.g. \cite{BR1}.
Moreover, the measure $\m_{\a\b}$
is pseudo--supported on the extreme point 
$\partial K(\car_{\a\b},\a^{\f_{\a\b}})$. 
Take 
$$
\Om_{\a\b}:=\partial K(\car_{\a\b},\a^{\f_{\a\b}})\,,
$$ 

and define on $M\cap\Om_{\a\b}$, 
$$
\tilde\n_{\a\b}(M\cap\Om_{\a\b}):=\m_{\a\b}(M)
$$
where $M$ is a Baire measurable set of 
$K(\car_{\a\b},\a^{\f_{\a\b}})$. Let $\n_{\a\b}$ be the completion of 
the probability measure $\tilde\n_{\a\b}$. 

Consider the integrable field $\{\gar_{\ff}\}_{\ff\in\Om_{\a\b}}$ of 
left Hilbert algebras whose elements have the form
$$
\ff\in\Om_{\a\b}\mapsto\pi_{\ff}(T)\Psi_{\ff}\,,\quad 
T\in\car_{\a\b}\,,
$$
where $(\pi_{\ff},\ch_{\ff},\Psi_{\ff})$ is the GNS representation of 
the modular state $\ff$ on $\gar_{\a\b}$. 

Let $\{\gam_{\ff\circ\eps_{1}}\}_{\ff\in\Om_{\a\b}}$, 
$\{\gn_{\ff\circ\eps_{2}}\}_{\ff\in\Om_{\a\b}}$, 
$\{\ga_{\ff\circ\eps}\}_{\ff\in\Om_{\a\b}}$ be the integrable fields of left Hilbert 
algebras analogously obtained starting from the states 
$\ff\circ\eps_{1}\in K(\cam_{\a\b},\a^{\f_{\a\b}\circ\eps_{1}})$
$\ff\circ\eps_{2}\in K(\cn_{\a\b},\a^{\f_{\a\b}\circ\eps_{2}})$
$\ff\circ\eps\in K(\ca_{\a\b},\a^{\f_{\a\b}\circ\eps})$
respectively.\footnote{Notice that the pull--back measures 
$\eps_{1}^{*}(\n_{\a\b})$, $\eps_{2}^{*}(\n_{\a\b})$, 
$\eps^{*}(\n_{\a\b})$ are precisely the orthogonal measures on 
$K(\cam_{\a\b},\a^{\f_{\a\b}\circ\eps_{1}})$, 
$K(\cn_{\a\b},\a^{\f_{\a\b}\circ\eps_{2}})$,
$K(\ca_{\a\b},\a^{\f_{\a\b}})$ corresponding to the abelian algebra
$Z_{\a\b}\sim I_{\a}\otimes Z_{\a\b}\otimes I_{\b}$ which is a common subalgebra of 
$Z(M_{\a\b})\sim Z(M_{\a\b})\otimes I_{\b}$, 
$Z(N_{\a\b})\sim I_{\a}\otimes Z(M_{\a\b})$, $Z(A_{\a\b})$ 
respectively, see e.g. \cite{A, BR1}. They provide the subcentral 
disintegration of the KMS states $\f_{\a\b}\circ\eps_{1}$, 
$\f_{\a\b}\circ\eps_{2}$, 
$\f_{\a\b}\circ\eps$ respectively.} 

Define
\begin{equation}
\label{dsm}
\Om:=\stackrel{\circ}{\bigcup}_{\a\b}\Om_{\a\b}\,,\quad 
\n:=\sum_{\a\b}\n_{\a\b}\,, 
\end{equation}
where `` $\stackrel{\circ}{\bigcup}$ '' stands for the disjoint union, and 
$\n_{\a\b}$ is understood as a measure on all of $\Om$ with support 
$\Om_{\a\b}$.
\begin{prop}
\label{6a1}
We have for the tensor product,
$$
M\ots_{R}N=\int^{\oplus}_{\Om}M(\ff)\ots_{R(\ff)}N(\ff)\di\n(\ff)
$$
where $\{R(\ff)\}_{\ff\in\Om}$ is the measurable field relative to 
the factor decomposition of $R$, and $\{M(\ff)\}_{\ff\in\Om}$, 
$\{N(\ff)\}_{\ff\in\Om}$ to the subcentral decomposition of 
$M$, $N$  
w.r.t. ${\displaystyle\bigoplus_{\a\b}I_{\a}\otimes Z_{\a\b}\subset 
Z(M)}$, 
${\displaystyle\bigoplus_{\a\b}Z_{\a\b}\otimes I_{\b}\subset Z(N)}$
respectively.
\end{prop}
\begin{pf}
We first consider the homogeneous case. Taking into account 
\eqref{11a1} and Proposition \ref{regal}, we have the inclusions
\begin{align*}
&I_{1}\otimes\cl(\gar_{\f})\subset\cl(\gam_{\f\circ\eps_{1}})
\subset B(H_{1})\otimes\cl(\gar_{\f})\,,\\
&\cl(\gar_{\f})\otimes I_{2}\subset\cl(\gn_{\f\circ\eps_{2}})
\subset\cl(\gar_{\f})\otimes B(H_{2})\,,
\end{align*}
$$
M\ots_{R}N\equiv\cl(\ga_{\f\circ\eps})=
\big(\cl(\gam_{\f\circ\eps_{1}})\otimes I_{2}\big)\bigvee
\big(I_{1}\otimes\cl(\gn_{\f\circ\eps_{2}})\big)\,.
$$

Taking into account the above considerations together with the 
structure of the commutant of the left von Neumann algebra of a left 
Hilbert algebra (see Section 3), we obtain by applying 
the results of \cite{R} and the Commutator Theorem 5.9 of \cite{SZ},
\begin{align*}
M\ots_{R}N\equiv&\cl(\ga_{\f\circ\eps})
=\int^{\oplus}_{\Om}\cl(\ga_{\ff\circ\eps})\di\n(\ff)\\
=&\int^{\oplus}_{\Om}\big(\cl(\gam_{\ff\circ\eps_{1}})\otimes I_{2}\big)\bigvee
\big(I_{1}\otimes\cl(\gn_{\ff\circ\eps_{2}})\big)\di\n(\ff)\\
\equiv&\int^{\oplus}_{\Om}M(\ff)\ots_{R(\ff)}N(\ff)\di\n(\ff)\,.
\end{align*}

Here, $M(\ff):=\cl(\gam_{\ff\circ\eps_{1}})$,
$N(\ff):=\cl(\gn_{\ff\circ\eps_{2}})$, $R(\ff):=\cl(\gar_{\ff})$.
The proof follows summing up all the homogeneous components appearing 
in th standard representation of $M\ots_{R}N$.
\end{pf}

We are ready to describe the structure of $M\ots_{R}N$ and its 
predual in the $\s$--finite case. 
\begin{thm}
\label{4}
Let $\{M(\ff)\}_{\ff\in\Om}$, $\{N(\ff)\}_{\ff\in\Om}$, 
$\{R(\ff)\}_{\ff\in\Om}$ be the measurable 
fields of von Neumann algebras appearing in the decompositions 
of $M$, $N$, $R$ given in Proposition \ref{6a1}, respectively. Define 
$\{\widetilde{M}(\ff)\}_{\ff\in\Om}$ 
and $\{\widetilde{N}(\ff)\}_{\ff\in\Om}$ as the measurable 
fields of von Neumann algebras such that
\begin{align*}
&\widetilde{M}(\ff)\otimes I(\ff)=M(\ff)\bigwedge 
\big(B(H_{\a\b})\otimes I(\ff)\big)\, ,\\
&I(\ff)\otimes\widetilde{N}(\ff)=N(\ff)\bigwedge \big(I(\ff)\otimes 
B(H_{\a\b})\big)
\end{align*}
whenever $\ff\in\Om_{\a\b}$.

Then we have
$$
M\ots_{R}N=\int^{\oplus}_{\Om}\widetilde{M}(\ff)\ots R(\ff)\ots
\widetilde{N}(\ff)\di\n(\ff)\,.
$$
\end{thm}
\begin{pf}
Looking at each fiber, we have by Lemma 1 of \cite{F2},
$$
M(\ff)\ots_{R(\ff)}N(\ff)=\widetilde{M}(\ff)\ots R(\ff)\ots_{R(\ff)}
R(\ff)\ots\widetilde{N}(\ff)\,,
$$
where the last equality follows as all the $R(\ff)$ are factors.
The proof follows as $R(\ff)\ots_{R(\ff)}R(\ff)$ coincides with $R(\ff)$.
\end{pf}

As an immediate corollary we have the structure of the predual
$\big(M\ots_{R}N\big)_{*}$, taking into account the description given 
in Section 3, of 
the predual of a direct integral of left von Neumann algebras arising 
from an integrable field of left Hilbert algebras.
\begin{cor}
In the situation of Theorem \ref{4}, we have
\begin{equation}
\label{71}
\big(M\ots_{R}N\big)_{*}=\int^{\oplus}_{\Om}
\widetilde{M}(\ff)_{*}\widehat{\otimes} R(\ff)_{*}
\widehat{\otimes}\widetilde{N}(\ff)_{*}\di\n(\ff)
\end{equation}
where $\widehat{\otimes}$ denotes the operator--projective tensor product 
between operator spaces given in \cite{ER1}.
\end{cor}
\begin{pf}
The proof directly follows from Theorem \ref{4}, taking into account 
Theorem 3.2 of \cite{ER2}.
\end{pf}

\bigskip\bigskip


\bigskip

\begin{thebibliography}{99}

\bibitem{A} Alfsen E. M.
{\it Compact convex sets and boundary integral},
Springer, Berlin Heidelberg New York, 1971.

\bibitem{BR1} Bratteli O., Robinson D. W.
{\it Operator algebras and quantum statistical mechanics I},
Springer, Berlin Heidelberg New York, 1981.

\bibitem{BR2} Bratteli O., Robinson D. W.
{\it Operator algebras and quantum statistical mechanics II},
Springer, Berlin Heidelberg New York, 1981.

\bibitem{BDL} Buchholz D., D'Antoni C., Longo R.
{\it Nuclear maps and modular 
structures II}, Commun. Math. Phys. {\bf 129} (1990), 115--138.

\bibitem{E} Effros E. , Ruan Z.-J.  {\it Operator spaces}, London 
Math. Soc. Monographs, New Series {\bf 29}, Clarendon Press, 
Oxford University Press (2000).

\bibitem{ER1} Effros E., Ruan Z.-J.  
{\it A new approach to operator spaces}, Canadian Math. Bull. {\bf 34} 
(1991), 329--337.

\bibitem{ER2} Effros E., Ruan Z.-J. 
{\it On approximation properties for operator spaces}, 
Internat. J. Math. {\bf 1} (1990), 163--187.

\bibitem{F1} Fidaleo F.
{\it Operator space structures and the split property}, 
J. Operator Theory {\bf 31} (1994), 207--218.

\bibitem{F} Fidaleo F. 
{\it On the split property for inclusions of $W^{*}$--algebras}, 
Proc. Amer. Math. Soc. {\bf 130} (2002), 121--127.

\bibitem{F2} Fidaleo F. {\it The predual of $W^{*}$--tensor products 
over $W^{*}$--subalgebras (separable case)}, 
J. Funct. Anal. {\bf 209} (2004), 194--205.

\bibitem{J} Goodman F. M., de la Harpe P., Jones V. F. R. 
{\it Coxeter graphs and towers of algebras}, Springer, Berlin Heidelberg New 
York, 1989.

\bibitem{L} Lance C. {\it Direct integral of left 
Hilbert algebras}, Math. Ann. {\bf 216} (1975), 11--28.

\bibitem{R1} Riedel N. {\it Topological direct integral of left 
Hilbert algebras I, II}, J. Operator Theory {\bf 5} (1981), 29--45, 
{\it ibidem} 213--229.

\bibitem{R} Riedel N. {\it Disintegration of KMS states and reduction 
of standard von Neumann algebras}, Pacific J. Math. {\bf 111} (1984), 
415--431.

\bibitem{St} Str\v{a}til\v{a} S.
{\it Modular theory in operator algebras}, Abacus press,
Tunbridge Wells, Kent, (1981).

\bibitem{SZ1} Str\v{a}til\v{a} S., Zsid\'o, L.
{\it Lectures on von Neumann algebras}, Abacus press, Tunbridge
Wells, Kent, (1979).

\bibitem{SZ} Str\v{a}til\v{a} S., Zsid\'o, L.
{\it The commutation theorem for tensor product over von Neumann 
algebras}, J. Funct. Anal. {\bf 165} (1999), 293--246.

\bibitem{Su} Summers S. J. {\it On the independence of local algebras in
Quantum Field Theory}, Rev. Math. Phys. {\bf 2} (1990), 201--247.

\bibitem{T} Takesaki M. 
{\it Theory of operator algebras I}, Springer, Berlin Heidelberg New 
York, 1979.

\bibitem{VW} Vesterstr{\o}m J., Wils W. {\it Direct integral of  
Hilbert spaces II}, Math. Scand. {\bf 26} (1970), 89--102.

\bibitem{W} Wils W. {\it Direct integral of  
Hilbert spaces I}, Math. Scand. {\bf 26} (1970), 73--88.

\end{thebibliography}
\end{document}